\documentclass[aop,MSNbibl,dvips]{arximspdf}

\doi{10.1214/12-AOP772} 
\volume{41}
\issue{5}
\pubyear{2013}
\firstpage{3606}
\lastpage{3616}

\makeatletter

\newtheorem{theorem}{Theorem}[section]
\newtheorem{lemma}[theorem]{Lemma}
\newtheorem{proposition}[theorem]{Proposition}
\newtheorem{corollary}[theorem]{Corollary}

\newproclaim{remark}[theorem]{Remark}

\newcommand{\rE}{\mathbb{E}}
\newcommand{\Cov}{\operatorname{cov}}
\newcommand{\ri}{\mathrm{i}}
\newcommand{\rd}{\mathrm{d}}

\newcommand{\var}{\operatorname{var}}
\newcommand{\sgn}{\operatorname{sgn}}

\makeatother

\begin{document}
\begin{frontmatter}

\title{Variance of partial sums of stationary sequences}
\runtitle{Variance of partial sums}

\begin{aug}
\author[A]{\fnms{George} \snm{Deligiannidis}\ead[label=e1]{deligian@stats.ox.ac.uk}}
\and
\author[B]{\fnms{Sergey} \snm{Utev}\corref{}\ead[label=e2]{sergey.utev@nottingham.ac.uk}}
\runauthor{G. Deligiannidis and S. Utev}
\affiliation{University of Oxford and University of Nottingham}
\address[A]{Department of Statistics\\
University of Oxford\\
OX1 3TG\\
United Kingdom\\
\printead{e1}} 
\address[B]{School of Mathematical Sciences\\
University of Nottingham\\
NG7 2RD\\
United Kingdom\\
\printead{e2}}
\end{aug}

\received{\smonth{10} \syear{2011}}

%
\begin{abstract}
Let $X_1, X_2, \ldots$ be a centred sequence of weakly stationary random
variables with spectral measure $F$ and partial sums $S_n = X_1 +
\cdots+ X_n$. We show that $\var(S_n)$ is regularly varying of index
$\gamma$ at infinity, if and only if $G(x):= \int_{-x}^x F(\rd x)$ is
regularly varying of index $2-\gamma$ at the origin ($0<\gamma<2$).
\end{abstract}

%

\end{frontmatter}

\section{Introduction}
Let $X_1, X_2, \ldots$ be a sequence of centered weakly stationary
random variables with finite second moments and spectral measure~$F$,
such that $r_k:= \operatorname{cov}(X_0,X_k) = \int_{-\pi}^{\pi}
\mathrm{e}^{\ri t k} \,\rd F(t)$, where to simplify calculations we assume that $F$ is a symmetric
measure about the origin, and let $G(x) = \int_{-x}^x F(\rd x)$.
Denote by $S_n$ the sequence of partial sums $S_n = X_1 + \cdots+ X_n$.

The main result of the paper is the following.
%
%
\begin{theorem}\label{thmlimit}
For $\gamma\in(0,2)$, define $C(\gamma) =\Gamma(1+\gamma)\sin
(\frac
{\gamma\pi}{2})/[\pi(2-\gamma)]$. Let $L(x)$ be a positive function,
slowly varying at infinity. Then:
\begin{longlist}
\item$G(x) \sim C(\gamma) K_0 x^{2-\gamma} L(1/x)$ as $x \to0$ if and
only if
\item$\var(S_n) \sim K_0 n^{\gamma} L(n)$ as $n\to\infty$.
\end{longlist}
In particular, $\var(S_n)/n\to K_0$ if and only if $G(x)/x\to K_0/\pi$.
\end{theorem}
The rate of growth of the variance of the partial sums $S_n$ has
received considerable attention in the literature due to its key role in
the limit theory of stationary random
sequences; see Bradley~\cite{Bradley07}, Chapter 8, and
Samorodnitsky~\cite{Samor2006}, Chapter~5, for comprehensive reviews.

\textit{Asymptotically linear behavior} $\var(S_n) \sim K_0 n$. To prove
asymptotic normality,
a common restriction on the dependence structure is to assume that the
growth of $\var(S_n)$ is asymptotically linear (Merlev{\`e}de, Peligrad
and Utev~\cite{Utev05}).

For this particular situation, there exist several results which
guarantee the convergence of $\var(S_n)/n$ under sufficient conditions
given in terms of mixing coefficients, linear dependence coefficients
or in terms of the covariances where it is well known that if $\lim_n
\sum_{k=-n}^n \rE X_0 X_k$ exists, then
$\lim_n \var(S_n)/n$ also exists in $[0,\infty)$, and the two limits
are equal (\cite{Bradley07}, Chapter 5).

In terms of the spectral measure $F$, Ibragimov's~\cite{Ibrag62} result
states that when $F$ is absolutely continuous, a sufficient condition
is the continuity of the spectral
density $f$ at the origin, in which case $\var(S_n)/n \to2\pi
f(0)$---this follows from the following representation:
\[
\frac{\var(S_n)}{n} = \frac{1}{n}\int_{-\pi}^{\pi}
\frac{\sin^2
({nt}/{2})}{\sin^2({t}/{2})} f(t) \,\rd t
\]
and Fej\'er's theorem on the Ces\`aro summability of the Fourier series
of the spectral
density $f$.

The continuity of the spectral density $f$ at the origin is by no means
necessary, and in fact Hardy and Littlewood~\cite{Hardy24}, Theorem C,
in 1924 proved that a necessary and sufficient condition is the
convergence of
\[
\frac{1}{t}\int_{-t}^{t} f(s) \,\rd s \to
c_0 \in(0,\infty),%
\]
which has appeared before in a probabilistic context~\cite{Bryc95}.

\textit{General case} $\var(S_n) \sim K_0 n^{\gamma} L(n)$. Necessary
conditions usually require restrictive assumptions on the covariances
such as regular variation in the Zygmund sense.
Several sufficient conditions are stated either in terms of the
covariances or of the spectral density,
for example, $f(x)\sim|x|^{-\alpha} L(1/|x|)$, for $\alpha\in(0,1)$,
implies that $\var(S_n)\sim n^{1+\alpha}L(n)$; see~\cite{Samor2006}.

Even when $\gamma=1$, the asymptotically nonlinear behavior of the variance
has appeared often in the limit theorems for dependent variables,
such as under general mixing conditions (see, e.g., \cite
{Ibrag62,Merlev06,Utev05}) or specific models such as
random walk in random scenery; see~\cite{Kest79,Bolt89}.

The case $\gamma>1$ frequently occurs in \textit{long-range dependent}
time-series when the covariances are not summable or the spectral
density has an appropriate singularity at the origin
which often results in non-Gaussian limiting behavior
\cite{Rosenblatt79,Taqqu75,Giraitis98,Surgailis00,Samor2006}. The case
$\gamma<1$ occurs when
the spectral density vanishes at the origin in which case both
non-Gaussian~\cite{Rosenblatt79}, and Gaussian limits
\cite{Taqqu75}, have appeared in the literature.

\textit{Technique.} It is not clear whether the approach of Hardy and
Littlewood~\cite{Hardy24} (and Zygmund~\cite{Zygmund02}), can handle
the slowly varying function and the case $\gamma\neq1$. We suggest
an alternative technique, which is based on weak convergence and
Fourier analysis of tempered distributions, which allows\vadjust{\goodbreak} us to work
directly with spectral measures, without assuming absolute continuity.

\textit{Subsequences.} Subsequences $\var(S_{2^n})/2^n$ have often been
applied through the
use of dyadic induction and stationarity, in the context of mixing
conditions, martingale
approximations, central limit theorems and invariance principles; see
\cite{Utev05}.
The question whether convergence along a subsequence is enough to
guarantee convergence of the full sequence has been around for some
time now, and it was presented to us as a conjecture by M. Peligrad.
Although the answer is positive under extra conditions such as $\rho$-mixing,
the necessary and sufficient condition stated in Theorem~\ref{thmlimit}
allows us to construct a counterexample proving that convergence along
dyadic subsequences does not imply convergence over the full sequence.
%
%
\begin{proposition}\label{thmcounterexample}
There exists a stationary process such that
$\var(S_{2^r})/2^{r}$ converges, but the full sequence $\var(S_{n})/n$
does not.
\end{proposition}
The proofs of Theorem~\ref{thmlimit} and Proposition~\ref{thmcounterexample}
are given in the next section along with several auxiliary results
which are of independent interest.
\section{Proofs}
We start by proving three auxiliary lemmas. By $C$ we denote a generic
positive constant.

\textit{Auxiliary results.}
Our starting point is the following inequality.
%
%
\begin{lemma}\label{propn}
For any $A>0$,
\[
\frac{4}{\pi^2} n^2 G(1/n) \leq\var(S_n) \leq G(\pi) +
\frac{\pi^2}{4} n^2 G(A/n) + \pi^2 \int
_{A/n}^{\pi} \frac{G(y)}{y^3} \,\rd y.
\]
\end{lemma}
\begin{pf}
Define the positive Fej\'er kernel $I_n(y)= \sin^2(ny/2)/\sin^2(y/2)$.
To prove the lower bound, notice that $I_n(y)\geq4n^2/ \pi^2$ for
$0<y<1/n$, and hence
\[
\var(S_n) = \int_0^{\pi}
I_n(y) G(\rd y)\geq\int_0^{1/n}
\frac
{4}{\pi^2}n^2 G(\rd y) \geq\frac{4}{\pi^2} n^2
G(1/n).
\]
To prove the upper bound, let $A\leq n$ and apply the bounds
$I_n(y)\leq n^2\pi^2/4$ for $y\leq A/n$ and $I_n(y)\leq\pi^2/y^2$ for
$y\geq A/n$ and integration by parts, to derive
\begin{eqnarray*}\hspace*{55pt}
\var(S_n) &=& \int_0^{A/n}
I_n(y) G(\rd y) + \int_{A/n}^{\pi}
I_n(y) G(\rd y)
\\
&\leq&\int_0^{A/n} \frac{n^2\pi^2}{4} G(\rd y) +
\int_{A/n}^{\pi} \frac
{\pi^2}{y^2} G(\rd y)
\\
&\leq&\frac{\pi^2}{4} n^2 G(A/n) + G(\pi) + \pi^2 \int
_{A/n}^{\pi
} \frac
{G(y)}{y^3} \,\rd y.\hspace*{55pt}\qed
\end{eqnarray*}
\noqed\end{pf}
The next result establishes that upper bounds of $\var(S_n)/g(n)$,
where $g(n) = n^{\gamma} L(n)$ for $\gamma\in(0,2)$, are equivalent
to upper bounds for the spectral measure $G$.
%
%
\begin{lemma}\label{lemma}
Suppose $\{n_k\}_{k\geq0}$ is a positive nondecreasing integer
sequence such that $n_k \to\infty$, and $\sup n_{k+1}/n_k = \kappa<
\infty$. Then the following are equivalent:

\begin{longlist}[(2)]
\item[(1)]
$\exists C>0$ and such that $\var(S_{n_k}) \leq C g(n_k)$;

\item[(2)] $\exists C>0$ such that $G(x) \leq C x^{2-\gamma} L(1/x)$;

\item[(3)] $\exists C>0$ such that $\var(S_{n}) \leq C g(n)$.
\end{longlist}
\end{lemma}
\begin{pf}
(1)$\Rightarrow$(2).
From Lemma~\ref{propn} and our assumptions, we have that for some positive
constant $C>0$,
$G(1/n_k) \leq
\pi^2Cn_k^{\gamma-2} L(n_k)/4$.
Thus, by monotonicity of $G$
and properties of slowly varying functions~\cite{Bin87},
for $1/n_{k+1} < x \leq1/n_k$,
\[
G(x) \leq G(1/n_k) \leq\frac{\pi^2 C}{4} \kappa^{2-\gamma}
x^{2-\gamma} L(1/x)\sup_{1/\kappa\leq\lambda\leq1} \frac{L(\lambda
/x)}{L(1/x)} \leq
C' x^{2-\gamma} L(1/x).
\]

(2) $\Rightarrow$ (3).
We apply Lemma~\ref{propn} with $A=1$ to get
\begin{eqnarray*}
\var(S_n) &\leq& G(\pi) + \frac{\pi^2}{4} n^2 G(1/n) +
\pi^2 \int_{1/n}^{\pi
} \frac
{G(y)}{y^3}
\,\rd y
\\
&\leq& C \biggl( 1 + g(n) + \int_{1/n}^{\pi}
y^{-\gamma-1} L(1/y) \,\rd y \biggr).
\end{eqnarray*}
Using the change of variables $x=1/y$ and since $\gamma-1 > -1$,
\[
\int_{1/n}^\pi y^{-\gamma-1} L(1/y) \,\rd y =
\int_{1/\pi}^n x^{\gamma-1}L(x) \,\rd x\sim
\frac{n^{\gamma
}L(n)}{\gamma}
\]
as $n\to\infty$, by the Tauberian theorem (\cite{Bin87}, Proposition 1.5.8).
Therefore there is a constant $C$ such that $\var(S_n) \leq C g(n)$
which completes
the proof since (3) $\Rightarrow$ (1) is obvious.
\end{pf}
The situation is similar when one considers lower bounds for $\var
(S_n)$ and $G(x)$.
%
%
\begin{lemma}\label{thmlowerbound}
Suppose that there exist positive constants $C_1$ and $C_2$ such that
$G(x) \leq C_1 x^{2-\gamma}L(1/x)$ and $\var(S_n) \geq C_2 g(n)$. Then
there exists a positive constant $C_3$ such that $G(x) > C_3x^{2-\gamma
}L(1/x)$.
\end{lemma}
\begin{pf}
We proceed by contradiction and assume that there is a sequence $0<y_k\to0$
such that $G(y_k)/y_k^{2-\gamma} L(1/y_k) \to0$, as $k \to\infty$.
Then we can construct a further sequence $1\leq A_k \to\infty$ slowly
enough so that
$A_k y_k \to0$, $G(A_k y_k)/y_k^{2-\gamma} L(1/y_k) \to0$, and
$L(1/A_ky_k)/L(1/y_k)\to1$.\vadjust{\goodbreak}

Then for $n_k = [1/y_k] +1\to\infty$, we have by Lemma~\ref{propn},
the Tauberian theorem
and the monotonicity of $G$, for generic positive constants $C, C'>0$,
as $k\to\infty$,
\begin{eqnarray*}
\frac{\var(S_{n_k})}{g(n_k)} &\leq& C \biggl( \frac{1}{g(n_k)} + \frac
{n_k^2 G(A_k/n_k)}{n_k^{\gamma} L(n_k)} +
\frac
{1}{n_k^{\gamma} L(n_k)}\int_{A_k/n_k}^{\pi}
{y_k}^{-3}G(y) \,\rd y \biggr)
\\
&\leq& C' \biggl(\frac{1}{g(n_k)} + \frac{G(A_ky_k)}{{y_k}^{2-\gamma
}L(1/y_k)} +
\frac{L(1/y_kA_k)}{A_k^{\gamma} L(1/y_k)} \biggr) \to0,
\end{eqnarray*}
which contradicts the assumptions of the lemma.
\end{pf}
%
%
\begin{remark}From Lemma~\ref{propn} it follows that the converse of
Lemma~\ref{thmlowerbound} is also true.
\end{remark}
%
%
\begin{remark}
For the boundary\vspace*{1pt} case $\gamma=0$, Theorem~\ref{thmlimit}
does not hold in general. For example, for $G(x)=x^2$ the direct
calculations show that
$\var(S_n)=4\ln(n)+O(1)$ which is not bounded. Actually,
Robinson \cite
{Robinson1960} proved that
\[
\sup_n\var(S_n)<\infty\quad\mbox{iff}\quad \int
_{0}^\pi x^{-2}\,\rd G(x)<\infty.
\]
When $\Cov(X_0,X_n)\to0$, the Leonov dichotomy holds; either $\var
(S_n)\to\infty$ or
$\sup_n\var(S_n)<\infty$; see Bradley~\cite{Bradley07}, Chapter 8.
However, the dichotomy is not true in general even for ergodic sequences
as it follows, for example, from the Aaronson and Weiss~\cite{AW2000}
construction on Chacon's ergodic transformations.

A nonergodic counterexample for Gaussian measures easily follows by
taking $G(\{2\pi2^{-k}\})=4^{-k}$
$k\geq2$, and following the calculations similarly to Samorodnitsky;
see~\cite{Samor2006}, Chapter 5. More exactly, then
\[
\sup_k\var(S_{2^k})<\infty\quad\mbox{and}\quad
\sup_n\var(S_{n})=\infty.
\]
\end{remark}
%
%
\begin{remark}
For the boundary case $\gamma=2$, the following
dichotomy easily follows from Lemma~\ref{propn}.
%
%
\begin{corollary}\label{ergodic}
Either $\liminf_{n\to\infty}\var(S_n)/n^2>0$ or $\var(S_n)/n^2 \to0$.
\end{corollary}
This fact also follows from the von Neumann $L_2$ ergodic theorem which
states that $\rE(S_n^2)/n^2$ vanishes
if and only if the spectral measure has no atom at the origin; see
\cite
{Doob90}.

Also from Lemma~\ref{propn}, the following corollary easily follows.
%
%
\begin{corollary}\label{another}
Let $\gamma\in(0,2]$, and $\{n_k\}_{k\in\mathbb{Z}}$ a nonnegative
increasing integer sequence such that $\sup_k n_{k+1}/n_k < \infty$.
Then the following are equivalent:
\begin{longlist}[(2)]
\item[(1)] $\var(S_{n_k})/n_k^\gamma\to0$;
\item[(2)] $x^{2-\gamma} G(x)\to0$ as $0<x \to0$; and
\item[(3)] $\var(S_{n})/n^\gamma\to0$.
\end{longlist}
\end{corollary}
Unlike the case $\gamma=2$, for $\gamma\neq2$, the equivalence does
not hold in general without the
assumption $\sup_k n_{k+1}/n_k < \infty$ as it follows from Theorem
\ref
{thmlimit} by using a slowly varying function $L$ such that $\liminf
_{n\to\infty}L(n)=0$ and
$\limsup_{n\to\infty}L(n)=\infty$.
\end{remark}

\textit{Proof of the main results.}
We are now ready to prove the main results.
Let $L$ be a positive function slowly varying at infinity,
$2-\gamma\in(0,2)$ and $g(n)=n^{\gamma}L(n)$.
The sufficiency part in Theorem~\ref{thmlimit} will be stated as an
independent lemma.
%
%
\begin{lemma}\label{fejerlemma}
Let $G(x) \sim x^{2-\gamma} L(1/x)$ as $x\to0$. Then, $\var
(S_{n})/g(n)\to1/C(\gamma)$.
\end{lemma}
\begin{pf}
Start with the representation
\begin{eqnarray*}
\var(S_{n}) 
&=& \int_{0}^{M}
\frac{\sin^2(y)}{n^2\sin^2(y/n)} n^2 G (2 \,\rd y/n) + \int_{M}^{n\pi/2}
\frac{\sin^2(y)}{n^2\sin^2(y/n)} n^2 G (2 \,\rd y/n)
\\
&=&\!: I_{n,M}+ J_{n,M}
\end{eqnarray*}
for fixed $M\leq n$.

The inequalities $n^2 \sin^2(y/n) \geq4y^2/\pi^2$, $G(x) \leq C
x^{2-\gamma}L(1/x)$ and integration by parts give
\[
J_{n,M} \leq\frac{\pi^2}{4} \int_{M}^{n\pi/2}
y^{-2} n^2 G (2 \,\rd y/n) \leq\frac{\pi^2 }{4} \biggl[ 4
\frac{G(\pi)}{\pi^2} + C 2^{3-\gamma} \int_{M}^{\infty}
n^{\gamma}\frac{L(n/2y)}{y^{1+\gamma}} \,\rd y \biggr].
\]
%
Bounding the integral term by using the change of variables $x=n/2y$
and the Tauberian theorem, we then derive
\[
J_{n,M} \leq\frac{\pi^2 }{4} \biggl[ 4\frac{G(\pi)}{\pi^2} + C
2^{3-\gamma} \frac{2^{\gamma}}{\gamma} \biggl(\frac{n}{2M}
\biggr)^{\gamma} L \biggl( \frac
{n}{2M} \biggr) \biggr] =O(1) + O
\bigl( g(n/2M)\bigr).
\]
By regular variation
$g(n/2M)/g(n) \to(2M)^{-\gamma}$ and therefore
%
%
\begin{equation}
\label{oct5eq1} \frac{\var(S_n)}{g(n)} = \frac{I_{n,M}}{g(n)} + O\bigl
(1/g(n)\bigr) + O
\bigl(M^{-\gamma}\bigr).
\end{equation}
Notice that for $y\leq M \leq n$,
$\sin^2(y)/n^2\sin^2(y/n) =\sin^2(y)/y^2 + O(M^2/n^2)$,
and thus
\[
\frac{I_{n,M}}{g(n)} = \int_0^M
\frac{\sin^2(y)}{y^2} \frac
{n^{2-\gamma
} G(2 \,\rd y /n)}{L(n)} + O\bigl(M^2/n^{\gamma}
\bigr).
\]
By regular variation of $G$,
it follows that for $y\leq M$,
\[
\mu_n \bigl([0,y)\bigr):= \frac{n^{2-\gamma} G(2y/n)}{L(n)
(2M)^{2-\gamma}} \to\biggl(
\frac{y}{M} \biggr)^{2-\gamma},
\]
which defines a probability measure on $[0,M]$, and hence, by weak
convergence, since $\sin^2(y)/y^2$ is continuous and bounded, there
exists a sequence $\mathcal{E}_M (n) \to0$, as $n\to\infty$ for all
$M$, such that
\begin{eqnarray*}
\frac{I_{n,M}}{g(n)} &=& 2^{2-\gamma}(2-\gamma)\int_0^M
\frac{\sin^2(y)}{y^{1+\gamma}} \,\rd y + \mathcal{E}_M (n) + O\bigl
(M^{-\gamma}
\bigr)
\\
&=& 2^{2-\gamma}(2-\gamma)\int_0^\infty
\frac{\sin^2(y)}{y^{1+\gamma}} \,\rd y + \mathcal{E}_M (n) + O\bigl
(M^{-\gamma}
\bigr)
\\
&=& \bigl(1/C(\gamma)\bigr) + \mathcal{E}_M (n) + O
\bigl(M^{-\gamma}\bigr);
\end{eqnarray*}
see~\cite{Gelfand64} for the integral. This together with (\ref
{oct5eq1}) implies the lemma.
\end{pf}
%
\begin{pf*}{Proof of Theorem~\ref{thmlimit}}
Implication (2) $\Rightarrow$ (1) immediately follows from Lemma \ref
{fejerlemma}.

For (1) $\Rightarrow$ (2), let $t_j\to\infty$ be a positive
increasing integer sequence.
Similar to Lemma~\ref{fejerlemma}, we derive
\[
\frac{\var(S_{t_j})}{g(t_j)} = \int_0^M
\frac{\sin^2(y)}{y^2} \frac{t_j^{2-\gamma} G(2\,\rd
y/t_j)}{L(t_j)} + O\bigl(M^2/t_j^\gamma
\bigr) + O\bigl(M^{-\gamma}\bigr).
\]
For $y\leq M$ we have
\[
\frac{t_j^{2-\gamma} G(2M/t_j)}{L(t_j)} \leq C M^{2-\gamma} \frac
{L(t_j/2M)}{L(t_j)} \leq C
M^{2-\gamma}.
\]
Helly's principle and a diagonal argument imply that there exists a
monotone increasing function $h$, defined on $[0,\infty)$, and a
subsequence $j'$ such that
%
%
\begin{equation}
\label{Helly} F_{t_{j'}}(y):= \frac{t_{j'}^{2-\gamma} G(2y
/t_{j'})}{L(t_{j'})} \to h(y)
\end{equation}
as $j' \to\infty$ for all continuity points $y$ of $h$. Since $h(y)
\leq CM^{2-\gamma}$ for $y\leq M$, and $\sin^2(y)/y^2$ is continuous
and bounded on $[0,M]$, by weak convergence we have that
\[
\int_0^M \frac{\sin^2(y)}{y^2} F_{t_{j'}}(\rd y) \to\int_0^M
\frac{\sin^2(y)}{y^2} h(\rd y).
\]
Therefore, writing an identity for arbitrary $M>0$ and then letting
$M\to\infty$
\[
K_0 = \lim_{j'\to\infty} \frac{\var(S_{t_{j'}})}{g(t_{j'})} = \int
_0^M \frac{\sin^2(y)}{y^2}h(\rd y) + O
\bigl(M^{-\gamma}\bigr) = \int_0^{\infty}
\frac{\sin^2(y)}{y^2} h(\rd y).
\]
Let $[x]$ denote the integer part of $x$, and notice that from (\ref
{Helly}) and regular variation of $G$, we also have
\[
F\bigl([rt_{j'}]\bigr) = \frac{[rt_{j'}]^{2-\gamma}
G(2y/[rt_{j'}])}{L([rt_{j'}])} \to r^{2-\gamma}
h(y/r)
\]
as $j' \to\infty$ for arbitrary $r>0$ and all continuity points $y/r$
of $h$.
Since $\var(S_n)/g(n)$ converges on the full sequence, it follows that
\[
K_0 = \lim_{j'\to\infty} \frac{\var(S_{[rt_{j'}]})}{g([rt_{j'}])} = \int
_0^{\infty} \frac{\sin^2(y)}{y^2} r^{2-\gamma} h(\rd
y /r)
\]
for any $r>0$, implying that
%
%
\begin{equation}\label{fourier1}
\int_0^{\infty} \frac{\sin^2(rx)}{x^2} h(\rd x) =
r^{\gamma} K_0.
\end{equation}
For $y>0$, let
$\psi(y):= \lim_{N\to\infty} \int_y^N x^{-2} h(\rd x)$,
which is well defined since by integration by parts we have
\[
\int_y^{\infty} x^{-2} h(\rd x)
= \lim_{N\to\infty} \biggl[ 2 \int_y^N
\frac{h(x)-h(y)}{x^3} \,\rd x \biggr] = 2 \int_y^{\infty}
\frac{h(x)-h(y)}{x^3} \,\rd x < \infty.
\]
The idea is to identify $\psi$ from its Fourier transform using the
convolution-type equation (\ref{fourier1}); along these lines we
continue by calculating the sine-transform of $\psi$ by interchanging
the integrals, which is allowed since the positive function $\psi$ is
bounded away from $0$ and $|\psi(y)| \leq C y^{-\gamma}$,
\begin{eqnarray*}
\int_0^a \sin(ry) \psi(y) \,\rd y &=& \int
_0^a \biggl(\int_{x=y}^{\infty}
\sin(ry) \frac{h(\rd
x)}{x^2} \biggr) \,\rd y
\\
&=& \frac{2}{r}\int_{0}^{a}
\frac{\sin^2(rx/2)}{x^2} h(\rd x) + \frac
{2}{r}\int_a^{\infty}
\frac{\sin^2(ax/2)}{x^2} h(\rd x)
\\
&\to&\biggl(\frac{r}{2} \biggr)^{\gamma-1} K_0 =
\lim_{a\to\infty} \int_0^a \sin(ry) \psi(y)
\,\rd y
\end{eqnarray*}
as $a \to\infty$ for all $r>0$.

For $y<0$, we define $\psi(y) = - \psi(-y)$
so that for any $r\in\mathbb{R}$, we have
%
%
\begin{equation}
\label{fourier2} \lim_{a\to\infty}\int_{-a}^a
\sin(ry) \psi(y) \,\rd y = 2^{2-\gamma
}\sgn(r) |r|^{\gamma-1}
K_0.
\end{equation}
To identify the function $\psi$ and therefore $h$, we apply Fourier analysis
and treat $\psi$ as a distribution acting on the Schwartz space of test
functions $\mathcal{S} = \mathcal{S}(\mathbb{R})$ such that
$\sup_x | x^{\alpha} \phi_n^{(\beta)}(x)| <\infty$ for all
nonnegative integers $\alpha$, $\beta$.
More exactly, we define a linear functional on $\mathcal{S}$ by
\[
\Psi[\phi] = \int_{0}^{\infty} \psi(y) \bigl(\phi(y)
- \phi(-y)\bigr)\,\rd y,
\]
which is continuous since
\begin{eqnarray*}
\bigl|\Psi[\phi]\bigr| &\leq& 4\sup_y \bigl|\phi'(y)\bigr| \biggl(\int
_{0}^{1} \psi(y)y\,\rd y \biggr) + 4
\sup_y \bigl|y\phi(y)\bigr| \biggl(\int_{1}^{\infty}
\bigl(\psi(y)/y\bigr)\,\rd y \biggr)
\\
&\leq& C_\psi\Bigl(\sup_y \bigl|\phi'(y)\bigr| +
\sup_y \bigl|y\phi(y)\bigr|\Bigr).
\end{eqnarray*}
The next step is to calculate the Fourier transform of the tempered
distribution $\Psi$ through the formula
$\hat{\Psi}[\phi] = \Psi[\hat{\phi}]$.
Then given $\phi\in\mathcal{S}$ we have
\begin{eqnarray*}
\Psi[\hat{\phi}] &=& \int_{y=0}^\infty\psi(y) \biggl(
\int_{t=-\infty}^{\infty} \bigl(\mathrm{e}^{\ri
ty}-
\mathrm{e}^{-\ri ty}\bigr) \phi(t) \,\rd t \biggr)\,\rd y
\\
&=& \ri\int_{y=-\infty}^\infty\psi(y) \biggl(\int
_{t=-\infty
}^{\infty} \sin(yt) \phi(t) \,\rd t \biggr) \,\rd y.
\end{eqnarray*}
Observe that $|\sin(yt)|\leq|yt|$, $|\psi(y)|\leq C|y|^{-\gamma}$, and
$1-\gamma>-1$ and so
for fixed~$a$,
\[
\int_{-a}^a \int_{-\infty}^{\infty}
\bigl|\psi(y)\bigr| \bigl|\sin(yt) \bigr| \bigl|\phi(t)\bigr| \,\rd t \,\rd y \leq\int_{-a}^a
\int_{-\infty}^{\infty} |y|^{1-\gamma} \bigl|t\phi(t)\bigr| \,\rd t
\,\rd y < \infty.
\]
Therefore by Fubini's theorem,
\begin{eqnarray*}
\Psi[\hat{\phi}] &=& \ri\lim_{a\to\infty} \int_{-a}^a
\psi(y) \int_{t=-\infty
}^{\infty
}\sin(yt) \phi(t) \,\rd t \,\rd y
\\
&=& \ri\lim_{a\to\infty} \int_{-\infty}^{\infty}\phi(t)
\int_{y=-a}^a \psi(y) \sin(yt) \,\rd y \,\rd t.
\end{eqnarray*}
We next bound the integrand in order to use dominated convergence. Let
$\tau= [ta/\pi]$, and write
\begin{eqnarray*}
I :\!&=& \int_{y=0}^a \psi(y) \sin(yt) \,\rd y =
\frac{1}{t}\int_{x=0}^{ta} \sin(x) \psi(x/t)
\,\rd x
\\
&=& \frac{1}{t}\sum_{j=0}^{\tau-1} \int
_{j\pi}^{(j+1)\pi} \sin(x) \psi(x/t) \,\rd x +
\frac{1}{t}\int_{\tau\pi}^{ta} \sin(x) \psi(x/t)
\,\rd x.
\end{eqnarray*}
Since $\tau\pi$ is the largest multiple of $\pi$ less than $ta$, for
$x\in[\tau\pi, ta]$ $\sin(x)$ does not change sign and therefore
\[
\biggl| \frac{1}{t}\int_{\tau\pi}^{ta} \sin(x)
\psi(x/t) \,\rd x \biggr| \leq|t|^{\gamma- 1} \int_{\tau\pi}^{ta}
\frac{|\sin(x)| }{|x|^{\gamma}}\,\rd x \leq C |t|^{\gamma-1}.
\]
The other term can be written as an alternating sum
\[
Q:=\frac{1}{t} \sum_{j=0}^{\tau-1}(-1)^j
\int_{j\pi}^{(j+1)\pi} \bigl|\sin(x)\bigr| \psi(x/t) \,\rd x.
\]
From the fact that $\psi$ is decreasing, we can then show that
\begin{eqnarray*}
c_j:\!&=&\int_{j\pi}^{(j+1)\pi} \bigl|\sin(x)\bigr|
\psi(x/t) \,\rd x \geq\int_{j\pi}^{(j+1)\pi} \bigl|\sin(x)\bigr| \psi
\bigl((j+1)\pi/t \bigr) \,\rd x
\\
&=&\int_{(j+1)\pi}^{(j+2)\pi} \bigl|\sin(x)\bigr| \psi\bigl((j+1)\pi/t
\bigr) \,\rd x \geq\int_{(j+1)\pi}^{(j+2)\pi} \bigl|\sin(x)\bigr| \psi(x/t
) \,\rd x= c_{j+1},
\end{eqnarray*}
and thus the sum $Q$ is conditionally convergent and in absolute value
less than its first term,
\begin{eqnarray*}
|Q| &\leq& \frac{1}{t} \biggl|\int_{0}^{\pi} \sin(x)
\psi(x/t) \,\rd x \biggr| =\frac{1}{t} \int_{0}^{\pi}
\sin(x) \psi(x/t) \,\rd x
\\
&\leq& \frac{C}{t} \int_{0}^{\pi} \sin(x)
\frac{t^\gamma}{x^\gamma
} \,\rd x =C t^{\gamma-1}.
\end{eqnarray*}
Overall the above calculations imply that for all $a>0$,
\[
\biggl| \phi(t) \int_{y=-a}^a \sin(yt) \psi(y) \,\rd y \biggr|
\leq C \bigl|\phi(t)\bigr| t^{\gamma-1},
\]
where $C$ does not depend on $a$. Furthermore
since $\gamma-1 >-1$, the function $|\phi(t)| t^{\gamma-1}$ has at most
an integrable singularity at the origin and is integrable.
Therefore by dominated convergence and (\ref{fourier2}),
\begin{eqnarray*}
\hat{\Psi}[\phi] &=& \ri\int_{-\infty}^{\infty}\phi(t)
\biggl(\lim_{a\to\infty} \int_{y=-a}^a \psi(y)
\sin(yt) \,\rd y \biggr) \,\rd t
\\
&=& \int_{-\infty}^{\infty}\phi(t) \bigl(
\ri2^{2-\gamma} K_0 \sgn(t) |t|^{\gamma-1} \bigr) \,\rd t.
\end{eqnarray*}
Then inverting the Fourier transform of the distribution $\Psi$ (see
\cite{Gelfand64}, e.g.) we identify the function $\psi(y)$, and by
standard calculations $h(x)$, for $x,y>0$, $\gamma\in(0,2)$
\[
\psi(y)= K_0 D(\gamma) y^{-\gamma},\qquad h(x) = \bigl(\gamma/(2-
\gamma)\bigr)K_0 D(\gamma)x^{2-\gamma},
\]
where $D(\gamma)=\Gamma(\gamma) 2^{2-\gamma} \sin(\gamma\pi
/2)/\pi$.
We have shown that for every integer sequence $t_j \to0$, there exists
a subsequence $t_{j'} \to\infty$ such that for $x,r>0$ and $\gamma
\in(0,2)$
\[
\frac{[rt_{j'}]^{2} G(x/[rt_{j'}])}{g([rt_{j'}])} \to r^{2-\gamma}
h(x/2r) = \bigl(\gamma/(2-\gamma)\bigr)
K_0 D(\gamma) (x/2)^{2-\gamma}.
\]
From this, by standard limiting arguments, we now deduce that
\[
\lim_{x\to0} \frac{G(x)}{x^{2-\gamma}L(1/x)} = \bigl(\gamma/(2-\gamma
)\bigr)
K_0 D(\gamma) (1/2)^{2-\gamma} = C(\gamma)K_0,
\]
which proves the theorem.\vadjust{\goodbreak}
\end{pf*}
%
\begin{pf*}{Proof of Proposition~\ref{thmcounterexample}}
The proof is through a counterexample. Let $G(x) = 2^{-k}$,
for $x\in(2^{-(k+1)},2^{-k}]$, for $k\geq1$. Then obviously\break
$\lim_{x\to0}G(x)/x$ does not exist, as different subsequences give
different limits.
Therefore by Theorem~\ref{thmlimit}
the limit of the full sequence $\var(S_n)/n$ cannot exist.

On the other hand, by direct calculation on the subsequence $2^r$,
\[
\frac{\var(S_{2^r})}{2^r} 
=\sum_{k=1}^\infty
\sin^2\bigl(2^{r-k-1}\bigr) 2^{k+1-r} \to\sum
_{k=0}^\infty\frac{\sin^2 (2^k)}{2^k} + \sum
_{k=1}^\infty2^k \sin^2
\bigl(2^{-k}\bigr) \in(0,\infty),
\]
completing the proof of the proposition.
\end{pf*}

\section*{Acknowledgments}

We would like to thank Professor M. Peligrad and Professor R. Bradley
for useful discussions.



\printaddresses

\end{document}